\documentclass[preprint]{elsarticle}
\usepackage{amscd, amssymb, amsthm}
\usepackage[all]{xy}
\usepackage{color}
\usepackage{url}

\begin{document}
\newcommand{\mono}[1]{%
\gdef\puA{#1}}
\newcommand{\puA}{}
\newcommand{\faculty}[1]{%
\gdef\puC{#1}}
\newcommand{\puC}{}
\newcommand{\facultad}[1]{%
\gdef\puD{#1}}
\newcommand{\puD}{}
\newcommand{\N}{\mathbb{N}}
\newcommand{\Z}{\mathbb{Z}}
\newcommand{\gl}{\mbox{\rm{gldim}}}
\newtheorem{teo}{Theorem}[section]
\newtheorem{prop}[teo]{Proposition}
\newtheorem{lema}[teo] {Lemma}
\newtheorem{ej}[teo]{Example}
\newtheorem{obs}[teo]{Remark}
\newtheorem{defi}[teo]{Definition}
\newtheorem{coro}[teo]{Corollary}
\newtheorem{nota}[teo]{Notation}
\def\fin{\mbox{\rm{findim}}}
\newcommand{\fidim}{\mbox{\rm{$\phi$dim}}}
\newcommand{\psidim}{\mbox{\rm{$\psi$dim}}}
\def\mod{\mbox{\rm{mod}}}
\def\add{\mbox{\rm{add}}}
\def\maxi{\mbox{\rm{max}}}
\def\pd{\mbox{\rm{pd}}}
\def\id{\mbox{\rm{id}}}
\def\gd{\mbox{\rm{gldim}}}
\def\mini{\mbox{\rm{min}}}
\def\ext{\mbox{\rm{Ext}}}
\def\enn{\hbox{\rm{End}}}
\def\hm{\hbox{\rm{Hom}}}
\def\Ker{\hbox{\rm{ker}}}

\newenvironment{note}{\noindent Notation: \rm}



\title{Igusa-Todorov functions for Artin algebras}

\author{Marcelo Lanzilotta}
\address{Universidad de la Rep\'ublica, Facultad de Ingenier\'ia -  Av. Julio Herrera y Reissig 565, Montevideo, Uruguay} 
\ead{marclan@fing.edu.uy}

\author{Gustavo Mata}
\address{Universidad de la Rep\'ublica, Facultad de Ingenier\'ia -  Av. Julio Herrera y Reissig 565, Montevideo, Uruguay} 
\ead{gmata@fing.edu.uy}

\begin{abstract}{
\noindent In this paper we study the behaviour of the Igusa-Todorov functions for Artin algebras $A$ with finite injective dimension, and Gorenstein algebras as a particular case. We show that the $\phi$-dimension and $\psi$-dimension are finite in both cases. Also we prove that monomial, gentle and cluster tilted algebras have finite $\phi$-dimension and finite $\psi$-dimension.}
\end{abstract}


\maketitle

\section{Introduction}

One of the most important conjecture in the Representation Theory of Artin algebras is the finitistic conjecture. It states that the $\sup\{\pd(M): M$ is a finitely generated module of finite projective dimension$\}$  is finite. As an attempt to prove the conjecture Igusa and Todorov defined in  \cite{kn:igusatodorov} two functions from the objects of $\mod A$ (the category of right finitely generated modules over an Artin algebra $A$) to the natural numbers, which generalizes the notion of projective dimension. They are known, nowadays, as the Igusa-Todorov functions, $\phi$ and $\psi$. One of its nicest features is that they are finite for each module, and allow us to define the $\phi$-dimension and the $\psi$-dimension of an algebra. These are new homological measures of the module category. In particular it holds that $\fin(A)\leq \fidim(A) \leq \psidim(A)\leq\gl(A)$ and they all coincide in the case of algebras with finite global dimension.\\
In the work \cite{kn:FLM} the authors presents relations of the $\phi$-dimension and the bifunctors $\mbox{\rm{ Ext}}_A( -, -)$ and $\mbox{\rm{Tor}}^A(-, -)$, that is they show how to compute the $\phi$-dimension of a module using these bifunctors. They also prove that the finiteness of the this dimension is invariant for derived equivalence.
Recently various works were dedicated to study and generalize the properties of  these functions, see for instance \cite{kn:marchuard}, \cite{kn:hlm1}, \cite{kn:xu}. In particular in \cite{kn:xu} the Igusa-Todorov functions were defined for the derived category of an Artin algebra.\\
A first calculation of the values of the $\phi$-dimension and $\psi$-dimension has been done for radical square zero algebras in \cite{kn: LMM}, a special case of $\Omega^n$-finite representation type algebras. In this work we prove that for a $\Omega^n$-finite representation type algebra (see Definition \ref{definfinite}) the $\phi$-dimension and $\psi$-dimension are both finite. In particular any finite dimensional monomial algebra has finite $\phi$-dimension (Corollary \ref{monomial algebras}). Also we study the Igusa-Todorov functions for Artin algebras $A$ such that $\id (A) < \infty$, and as particular case we get the result for Gorenstein algebras (observe that this algebras are not always of $\Omega^n$-finite representation type).\\
This paper is organized as follows: after the introduction and the preliminary section devoted to fixing the notation and recalling the basic facts needed in this work, we get section 3 about Igusa-Todorov functions and $ ^\bot A$ (see Definition \ref{perpA}). Using results given in \emph{Stable module theory} \cite{kn:Auslander-Bridger}, we get the following:\\ 

{\underline{\bf{Theorem A}}:} $\fidim( ^\bot A) = \psidim ( ^\bot A) = 0$.\\

As a consequence of the previous theorem, we obtain the following result:\\

{\underline{\bf Corollary B}}: Let $A$ be an Artin algebra such that $\id (A) = n < \infty$, then 

$$\fin (A) \leq \fidim (A) \leq \psidim (A) \leq n.$$

In section 4 we get applications to Gorenstein algebras, gentle algebras and cluster tilted algebras. The main theorem of this section is:\\

{\underline{\bf Theorem C}}: If $A$ is a Gorenstein algebra, then 

\begin{enumerate}
\item $\fin (A) = \fidim (A) = \psidim (A) < \infty.$

\item Let $m$ be the minimun  integer such that $A$ is a $m$-Gorenstein algebra, then: 
$$\fin (A) = \fidim (A) = \psidim (A) = m.$$
\end{enumerate}

Because an Artin  algebra $A$ is selfinjective if and only if $A$ is $0$-Gorenstein, we get as a consequence of the previous theorem the main result of \cite{kn:marchuard}:\\

{\underline{\bf Corollary D}}: \cite[Corollary 6]{kn:marchuard}
Let $A$ be an Artin  algebra. Then $A$ is selfinjective if and only if $\fidim(A) = 0$.

\section{Preliminaries}

\subsection{Gorenstein modules and $ ^\bot A$}

\begin{defi}\label{perpA}

Let $A$ be an Artin algebra. We will denote with $ ^\bot A$ the full subcategory of $\mod A$ whose objects verify $\ext_A^i(M, A) = 0$ for $i\geq 1$.

\end{defi}

\begin{obs}\label{A perp}

Let $A$ be an Artin algebra. The subcategory $^\bot A$ is closed by kernel of epimorphisms, is closed by extensions and in particular is closed by finite direct sums of modules. In addition every direct summand for a module in $^\bot A$ is also a module in $^\bot A$.

\end{obs}

The following notation is in order to simplify.

\begin{nota}
We will denote with $(\cdot )^{\ast}$ the functor $\hm_A( \cdot, A)$

\begin{itemize}
  \item $M^{\ast} = \hm_{A}( M, A)$ when $M$ is an $A$-module.
  \item $f^{\ast} = \hm_{A}( f,A)$ when $f$ is a morphism of $A$-modules.
\end{itemize}

\end{nota}

\begin{defi}

A finitely generated $A$-module $G$ is \textbf{Gorenstein projective} if there exist a exact sequence:
$$ \xymatrix{ \ldots \ar[r] & P_{-2} \ar[r]^{p_{-2}} & P_{-1} \ar[r]^{p_{-1}} & P_0  \ar[r]^{p_0}  & P_1 \ar[r]^{p_1}& P_2 \ar[r]^{p_2} & \ldots }$$
\noindent such that $G \cong \Ker p_0 $,  $P_i$ is a projective module for $i \in \mathbb{Z}$ and if the functor $\hm_A( \cdot , A)$ is applied then it is obtained the following exact sequence:
$$ \xymatrix{ \ldots \ar[r]& {P_{2}}^{\ast} \ar[r]^{{p_{1}}^{\ast}} & {P_{1}}^{\ast} \ar[r]^{{p_{0}}^{\ast}} & {P_{0}}^{\ast} \ar[r]^{{p_{-1}}^{\ast}}  & {P_{-1}}^{\ast}\ \ar[r]^{{p_{-2}}^{\ast}}&  {P_{-2}}^{\ast} \ar[r]^{{p_{-3}}^{\ast}}& \ldots }$$
\end{defi}

The next properties are well known (see \cite{kn:Zhang}):

\begin{obs} Given $A$ an Artin algebra then: 

\begin{enumerate}
  \item Every direct sum of Gorenstein projective modules is a Gorenstein projective module.

  \item Every direct summand of a Gorenstein projective module is a Gorenstein projective module.

  \item Every projective module is a Gorenstein projective module.
  
  \item Every Gorenstein projective module is a projective module or its projective dimension is infinite. 

\end{enumerate}

\end{obs}

\begin{nota}

Let $A$ be an Artin algebra. The full subcategory of $\mod A$ which objects are Gorenstein projectives will be denoted by $\mathcal{G} P_A$. 
  
\end{nota}

For a proof of the next proposition see \cite{kn:Zhang}.

\begin{prop}\label{Gorenstein sub Aper} \cite{kn:Zhang}

$\mathcal{G} P_A $ is a full subcategory of $ ^\bot A$.\\

\end{prop}

\subsection{Igusa-Todorov functions}

In this section, we show some general facts about the Igusa-Todorov functions for an Artin algebra $A$. Our objective here is to introduce some material which we will use in the following sections.

\begin{defi}
Let $K_0$ denote the quotient of the free abelian group generated by one symbol  $[M]$,  for each isomorphism class of right finitely generated $A$-module, and relations given by:
\begin{enumerate}
  \item $[M]-[M']-[M'']$ if  $M \cong M' \oplus M''$.
  \item $[P]$ for each projective module.
\end{enumerate}
\end{defi}

\noindent  We use frequently the following observation.

\begin{obs}\label{positivecone}
It is clear that $K_0$ has as a basis the set of symbols $[M]$, one for each class of isomorphism of indecomposable non projective module. Moreover every element in $K_0$ can be written in the form $[M]-[N]$, for some pair of, not necessarly indecomposable, modules $M$ and $N$.
\end{obs}

\noindent Let $\overline{\Omega}: K_0 \rightarrow K_0$ be the group endomorphism induced by $\Omega$, and let   $K_i = \overline{\Omega}(K_{i-1})= \ldots = \overline{\Omega}^{i}(K_{0})$ for $i \in \mathbb{N}$. Now if $M$ is a finitely generated $A$ module then $\langle add M\rangle$ denotes the subgroup of $K_0$ generated by the classes of indecomposable summands of $M$.

\begin{defi}\label{monomorfismo}
The \textbf{(right) Igusa-Todorov function $\phi$} of $M\in \mod A$  is defined as  $\phi_{r}(M) = min\{l:
\overline{\Omega}{|}_{{\overline{\Omega}}^{l+s}\langle add M\rangle}$ is a monomorphism for all $s \in \mathbb{N}\}$.
\end{defi}

\noindent In an analogous way, using the cosyzygy, we may define the left Igusa-Todorov function $\phi_{l}(M)$. Using duality we can see that $\phi_{l}(M)=\phi_{r}(D(M))$, where $D$ is the usual duality between $\mod A$ and $\mod A^{op}$ (see \cite[II.3]{kn:ARS}). In case there is no possible misinterpretation we will use the notation $\phi$ for the right Igusa-Todorov function. We can find the next propositions in \cite{kn:igusatodorov} and \cite{kn:hlm1}.

\begin{prop}\cite{kn:igusatodorov}, \cite{kn:hlm1}\label{it1} \label{Huard1}

Let $M$ and $N$ be $A$-modules, then:

\begin{enumerate}
  \item $\phi(M) = \pd (M)$ if $M$ has finite projective dimension.
  \item $\phi(M) = 0$ if $M$ is indecomposable and has infinite projective dimension.
  \item $\phi(M) \leq \phi(M \oplus N)$.
  \item $\phi(M^{k}) = \phi(M)$ for $k \in \mathbb{N}$.
  \item If $M \in \mod A$, then $\phi(M) \leq \phi(\Omega(M))+1$.
\end{enumerate}

\end{prop}

\begin{defi}
The \textbf{(right) Igusa-Todorov function $\psi$} of $M\in \mod A$  is defined as  $\psi_{r}(M) = \phi_{r}(M) + \sup\{ \pd (N) : \Omega^{\phi (M)}(M) = N \oplus N'\mbox{ and } \pd (N)< \infty \}$.
\end{defi}

\noindent In an analogous way, using the cosyzygy, we may define the left Igusa-Todorov function $\psi_{l}(M)$. Using duality we can see that $\phi_{l}(M)=\phi_{r}(D(M))$, where $D$ is the usual duality between $\mod A$ and $\mod A^{op}$ (see \cite[II.3]{kn:ARS}). In case there is no possible misinterpretation we will use the notation $\psi$ for the right Igusa-Todorov function.

\begin{prop}\cite{kn:igusatodorov}, \cite{kn:hlm1} \label{it2} \label{Huard2}

Let $M$ and $N$ be $A$-modules, then:

\begin{enumerate}
  \item $\psi(M) = \pd (M)$ if $M$ has finite projective dimension.
  \item $\psi(M) = 0$ if $M$ is indecomposable and has infinite projective dimension.
  \item $\psi(M) \leq \psi(M \oplus N)$.
  \item $\psi(M^{k}) = \psi(M)$ for $k \in \mathbb{N}$.
  \item If $N$ is direct summand of $\Omega^{n}(M)$ where $n\leq \phi (M)$ and $\pd (N) < \infty$, then $\pd(N) + n \leq \psi (M)$.
  \item If $M \in \mod A$, then $\psi(M) \leq \psi(\Omega(M))+1$.
\end{enumerate}

\end{prop}

\begin{defi}

Let $A$ be an Artin algebra and $\mathcal{C}$ a full subcategory of $\mod A$. We define:

\begin{itemize}

\item $\fidim (A) = \sup \{\phi(M) $ such that $ M \in \mod A \}$.

\item $\psidim (A) = \sup \{\psi(M) $ such that $ M \in \mod A \}$.

\item $\fidim (\mathcal{C}) = \sup \{\phi(M) $ such that $ M \in Ob \mathcal{C} \}$.

\item $\psidim (\mathcal{C}) = \sup \{\psi(M) $ such that $ M \in Ob \mathcal{C} \}$.

\end{itemize}

\end{defi} 

\begin{obs} \label{inecuaciones} For an Artin algebra $A$ we have the following:
\begin{itemize}

\item $\fin (A) \leq \fidim (A) \leq \psidim (A)$.

\item $\fidim(A) < \infty$ if and only if $\psidim (A) < \infty$.

\end{itemize}
 
\end{obs}

\begin{defi}

Let $M$ be an $A$-module. We call it {\bf{$\Omega$-invariant module}} if $\Omega(M) \in \add M$.

\end{defi}

For a proof of the following proposition see \cite{kn: LMM}. 

\begin{prop} \cite[Proposition 3.6]{kn: LMM}\label{invariante}

Let $M =\oplus_{i = 1}^k M_i$ be an $\Omega$-invariant module where $M_i$ is indecomposable for $i = 1, \ldots , k$ with $[M_i] \neq [M_j]$ whenever $i \neq j$, then $\phi (M) \leq k$.

\end{prop}

\section{Igusa-Todorov functions and $ ^\bot A$}

\subsection{$\Omega^n$ finite representation type algebras}

\begin{defi}{\label{definfinite}}

Let $A$ be an Artin algebra. For $n \in \mathbb{N} $, we say that $K_n$ is of {\bf{finite representation type}} if there are only finitely many indecomposables modules $M$, up to isomorphism, such that $M \in \add N$ with $[N] \in K_n$. In this case $A$ will be called of {\bf $\Omega^n$-finite representation type}. 

\end{defi}

\begin{teo} \label{finito}

Let $A$ be an $\Omega^n$-finite representation type algebra, then $\fidim(A) < \infty$ and $\psidim(A) < \infty$.

\begin{proof}

By the previous hypothesis we can fix $\bar{M}$ an $A$-module, such that $\add N \subset \add \bar{M}$ for every $[N] \in K_n$.\\
Let $M$ be any $A$-module, then $\Omega^n(M) \in \add \bar{M}$ so by Proposition \ref{Huard1} $\phi(M) \leq \phi(\Omega^n(M))+n \leq \phi(\bar{M})+n$. The finiteness of the $\psidim (A)$ follows from Remark \ref{inecuaciones}.\end{proof}

\end{teo}

\subsection{Examples of $\Omega^n$-finite representation type algebras}

\begin{defi}

An $A$-module $M$ is said to be {\bf{uniserial}} if it has a
unique composition series. An algebra $A$ is said to be {\bf right serial} if every indecomposable projective right $A$-module is uniserial.

\end{defi}

The following theorem can be found in \cite{kn: Zimmermann2}.

\begin{teo} \cite[Proposition 1]{kn: Zimmermann2}\label{serial}
If $A$ is a right serial algebra, then every submodule $M$ of a projective module is a direct sum of uniserial right ideals. In particular $A$ is of $\Omega^1$-finite representation type.

\end{teo}

\begin{defi}

A $\Bbbk$-algebra is \textbf{special biserial} if it is morita equivalent to an algebra $\frac{\Bbbk Q}{I}$ where $I$ is an admissible ideal with the following conditions:

\begin{itemize}
  \item At each vertex at most two arrows begin.
  \item At each vertex at most two arrows end.
  \item For each arrow $\beta \in Q_1$ there is at most one arrow $\gamma \in Q_1$ such that $\beta \gamma \notin I$.
  \item For each arrow $\beta \in Q_1$ there is at most one arrow $\alpha \in Q_1$ such that $\alpha \beta \notin I$.
\end{itemize}

\end{defi} 
 
The following theorem can be found in \cite{kn:Ringel}.

\begin{teo}\cite[Section 5.11]{kn:Ringel}\label{sbiserial}

Let $A$ be a special biserial algebra without a projective
injective module, then A is $\Omega^1$-finite representation type.

\end{teo}  
 
\begin{coro}

If $A$ is right serial or special biserial algebra, then $\fidim (A) < \infty$ and $\psidim (A) < \infty$.

\begin{proof}

By Theorems \ref{serial} and \ref{sbiserial} any right serial or special biserial algebra is $\Omega^1$-finite representation type. The result follows using Theorem \ref{finito}.\end{proof}

\end{coro}

We can apply Proposition \ref{invariante} to monomial algebras because of the following theorem (see \cite{kn:Zimmermann}) and obtain Corollary \ref{monomial algebras}.

\begin{teo}\cite[Theorem 3.I]{kn:Zimmermann}\label{Z}
Let $A$ be a finite dimensional monomial algebra. If $M$ is an $A$-module such that $[M] \in K_t$ with $t \geq 2$ , then it is a direct sum of right ideals of $A$ which are generated by paths with length greater or equal than $1$ ($\mathbb{P}^{\geq 1}$). In particular monomial algebras are of $\Omega^2$-finite representation type.
\end{teo}

\begin{coro}\label{monomial algebras}

Let $A = \frac{\Bbbk Q}{I}$ be a monomial algebra, then $\fidim (A) \leq {\dim }_{\Bbbk}A- n+2$, where $n = |Q_0|$.

\begin{proof}

Given $M$ an $A$-module, then $\Omega^{2}(M)$ is a direct sum of right ideals of $A$ which are generated by paths of $Q$ with length greater or equal than $1$. If $N = \oplus_{p \in \mathbb{P}^{\geq 1}} \langle p\rangle $, then is easy to see that $rk (\langle \add N \rangle) = {\dim }_{\Bbbk} (A) -n$. Now using the fact that $\Omega (N) \in \add N$ (see Theorem \ref{Z}), by Proposition \ref{invariante} then we have  $\phi(N) \leq {\dim }_{\Bbbk} (A) -n$.\\
In the other hand $\Omega^{2}(M)$ is a direct summand of $N^l$, then $\phi(\Omega^{2}(M)) \leq \phi(N) \leq {\dim }_{\Bbbk} (A) -n$. Finally by Proposition \ref{Huard1} can be deduced $\phi(M) \leq {\dim }_{\Bbbk}A- n+2 $.\end{proof}

\end{coro}

\subsection{Modules in $ ^\bot A$ and Igusa-Todorov functions}

In this section we will compute the $\phi$-dimension and the $\psi$-dimension for a family of algebras independently if $K_n$ is of finite or infinite representation type with $n \in \mathbb{N}$.

\begin{defi}

Given $\mathcal{C}$ a full subcategory of $\mod A$, we define the projective stable category $\underline{\mathcal{C}}$ in the following way:

\begin{itemize}
  \item The objects in $\underline{\mathcal{C}}$ are the same that the objects in $\mathcal{C}$.
  \item The morphisms set $\hm_{\underline{\mathcal{C}}}(M, N) = \frac{\hm_{\mathcal{C}}(M, N)}{P(M,N)}$, where $P(M,N)$ is the subset of morphisms from $M$ to $N$ such that factorize by a projective module.
  \item The composition is the induced from $\mathcal{C}$.
\end{itemize}

\end{defi}

The following two lemmas will be useful (see \cite{kn:FLM} for the first one and \cite{kn:Auslander-Bridger} for the second one).

\begin{lema}\cite[Proposition 3.1]{kn:FLM}\label{equivalencia FLM}

Let $A$ be an Artin algebra, $M$, $N \in \mod A$, then the following are equivalent:

\begin{enumerate}
  \item $M \oplus P_0(N) \cong N \oplus P_0 (M)$ in $\mod A$.
  \item $M \cong N$ in $\underline{\mod} A $.
  \item $[M] = [N]$ in $K_0(A)$.
\end{enumerate}

\end{lema}

As we have seen in Remark \ref{A perp}, the subcategory $^\bot A$ is invariant by $\Omega$, therefore this induces a functor $\Omega : \underline{ ^\bot A} \rightarrow \underline{ ^\bot A}$.

\begin{lema} \cite{kn:Auslander-Bridger} \label{Auslander-Bridger}

The functor $\Omega : \underline{ ^\bot A} \rightarrow \underline{ ^\bot A}$ is faithfull and full.\\

\end{lema}

\begin{nota}
Let $\mathcal{C}$ be a subcategory of $\mod A$. The subgroup of $K_0$ formed by the classes of modules in $\mathcal{C}$, will be denoted by $[\mathcal{C}]$.  
\end{nota}

\begin{prop}\label{sizigia n-esima di A<n}

Let $A$ be an Artin algebra such that $\id A_A \leq n$, then $K_n \subset \  [ ^\bot \! A]$.

\begin{proof}

Let consider the following exact sequence:

$$ \xymatrix{ 0 \ar[r] & K \ar[r] & P_{n-1} \ar[r]^{p_{n-1}} &  \ldots \ar[r] & P_1  \ar[r]^{p_1}  & P_0 \ar[r]^{p_0}& M \ar[r] & 0}$$

\noindent where each $P_i$ is projective.
It follows that $\ext_A^{i}(K,A) = \ext_A^{i}(\Omega^{n}(M), A)\cong \ext_A^{n+i}(M, A)$. Now using the fact that $\id A \leq n$, it is clear that $\ext_A^{i}(K,A) = 0$ for $i \geq 1$.\end{proof}

\end{prop}

\begin{teo}\label{phi en perpA}

$\fidim( ^\bot A) = \psidim ( ^\bot A) = 0$.

\begin{proof}
Let us prove that the morphism $\overline{\Omega}{|}_{[ ^\bot A]}$ is injective and then, by Definition \ref{monomorfismo} and Remark \ref{A perp}, we obtain the result.\\
Suppose that $\ker \ \overline{\Omega}{|}_{[ ^\bot A]} \neq 0$, that means there are elements $[M_1]$ and $[M_2]$ in $[ ^\bot A]$ such that $\overline{\Omega}([M_1]-[M_2]) = 0$. Then $\overline{\Omega}([M_1]) = \overline{\Omega}([M_2])$ and this implies that $[\Omega(M_1)] = [\Omega(M_2)]$. By Lemma \ref{equivalencia FLM} it is equivalent to $\Omega(M_1) \cong \Omega(M_2)$ in $\underline{ ^\bot A}$.\\
Now, by Lemma \ref{Auslander-Bridger}, $\Omega : \underline{ ^\bot A} \rightarrow \underline{ ^\bot A}$ is a full and faithfull functor, so $M_1 \cong M_2$ in $\underline{ ^\bot A}$. Using Lemma \ref{equivalencia FLM} we get that $[M_1] = [M_2]$.\\
Because $\phi(M) = 0$, $M$ has no summands of finite projective dimension, then $\psi(M) = 0$.\end{proof}

\end{teo}

The following example shows that the inclusion $ ^\bot A \subset \{M \in \mod A$ such that $\phi (M) = 0\}$ can be strict.\\

\begin{ej}

Consider the following radical square zero algebra $A = \frac{\Bbbk Q}{J^{2}}$ where $Q$ is definded as follows:

$$
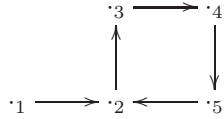
\xymatrix{ & \cdot_3 \ar[r] & \cdot_4 \ar[d] \\ \cdot_1  \ar[r] & \cdot_2 \ar[u] & \cdot_5 \ar[l] }$$

It is clear that $\ext_A^{1}(S_1, P_5) \neq 0$ because there exists the non split short exact sequence:

$$\xymatrix{0\ar[r] & P_5 \ar[r] & I_2 \ar[r] & S_1 \ar[r]& 0}$$

\noindent It is also clear that $\pd S_1 = \infty$, and follows that $\phi (S_1) = 0$. As a consequence we obtain $ ^\bot A \subsetneq \{ M \in \mod A$ such that $ \phi(M) = 0\}$.

\end{ej}

\begin{coro}\label{phidim diA<n}

Let $A$ be an Artin algebra such that $\id (A) = n < \infty$ then 

$$\fin (A) \leq \fidim (A) \leq \psidim (A) \leq n.$$

\begin{proof}

Given $M \in \mod A$, we know by Proposition \ref{sizigia n-esima di A<n} that $\Omega^n(M) \in ^\bot \!\!\! A$, then by Theorem \ref{phi en perpA}, $\phi(\Omega^n(M)) = \psi(\Omega^n(M)) = 0$. Therefore using Propositions \ref{Huard1} and \ref{Huard2} we obtain that $\phi(M) \leq \psi(M) \leq n$.\end{proof}

\end{coro}

\section{Applications to Gorenstein modules and Gorenstein algebras}
In this section we start proving as a corollary of Theorem \ref{phi en perpA} and Proposition \ref{Gorenstein sub Aper} that $\phi(G) = \psi(G) = 0$ for all $G$ Gorenstein projective module.
This result appeared in the list of abstracts of ICRA XV 2012, Bielefeld, but until now there is no published document with the proof of it.

\begin{coro}\label{Gorenstein proy}

$\fidim( \mathcal{G} P_A) = \psidim (\mathcal{G} P_A) = 0$.

\end{coro}

Our purpose now, is to compute the $\phi$-dimension and $\psi$-dimension of Gorenstein algebras. 

\begin{prop}\label{interseccion K_i}

Given $A$ an Artin algebra, we have that:

\begin{itemize}
\item $[\mathcal{G} P_A] \subseteq \cap_{i = 0}^{\infty} K_i$ and
\item $\overline{\Omega} ([\mathcal{G} P_A]) = [\mathcal{G} P_A]$.
\end{itemize}

\begin{proof}

Given $G$ an indecomposable Gorenstein projective non projective module, then $G \cong \Ker p_0 $ where $p_0$ is a morphism in the following exact sequence:

$$ \xymatrix{ \ldots \ar[r] & P_{-2} \ar[r]^{p_{-2}} & P_{-1} \ar[r]^{p_{-1}} & P_0  \ar[r]^{p_0}  & P_1 \ar[r]^{p_1}& P_2 \ar[r]^{p_2} & \ldots }$$

\noindent where $P_i$ is a projective module for all $i \in \mathbb{Z}$.

The next assertions follow:

\begin{itemize}
\item $[G] = \overline{\Omega} ([\Ker p_1])$, where $\Ker p_1$ is also Gorenstein projective, and
\item $\overline{\Omega}([G]) = [\Ker p_{-1}]$.
\end{itemize} 

Using the last facts, we obtain $\overline{\Omega} ([\mathcal{G} P_A]) = [\mathcal{G} P_A]$ and, in particular, $[\mathcal{G} P_A] \subseteq K_1$. We conclude that $[\mathcal{G} P_A] \subseteq \cap_{i = 0}^{\infty} K_i$.\end{proof}

\end{prop}

\begin{defi}

An Artin algebra $A$ is called $n$-Gorenstein if $\id(A_A) \leq n$ and $\pd( DA^{op}_A) \leq n$ with $n \in \mathbb{N}$. An Artin algebra $A$ is called Gorenstein if it is $n$-Gorenstein for some $n \in \mathbb{N}$.

\end{defi}

The following proposition can be seen in \cite{kn:Zhang}: 

\begin{prop} \cite[Corollary 3.4]{kn:Zhang} \label{Zhang} \label{sizigia n-esima Gor}

Let $A$ be a $n$-Gorenstein algebra, and 

$$ \xymatrix{ 0 \ar[r] & K \ar[r] & P_{n-1} \ar[r]^{p_{n-1}} &  \ldots \ar[r] & P_1  \ar[r]^{p_1}  & P_0 \ar[r]^{p_0}& M \ar[r] & 0}$$

\noindent be an exact sequence with $P_i$ projective, then $K$ is a Gorenstein projective $A$-module.

\end{prop}

As a direct consequence of the Propositions \ref{interseccion K_i} and \ref{sizigia n-esima Gor} we obtain the next result:

\begin{coro}

If $A$ is an $n$-Gorenstein algebra, then $[\mathcal{G} P_A] = K_n$.

\end{coro}

The following lemma shows that the concept of being n-Gorenstein is symmetric (see \cite{kn: happel}). 

\begin{lema} \cite[Lemma 1.2]{kn: happel} \label{Gorenstein reversible}

Let $A$ be a Gorenstein algebra and $\pd DA^{op} = r$, then $\id A = r$.\\

\end{lema}

The next result follows directly from Corollary \ref{phidim diA<n}.

\begin{teo}\label{phi algGorenstein}\label{psi algGorenstein}

If $A$ is a Gorenstein algebra, then:

\begin{enumerate}
\item $\fin (A) = \fidim (A) = \psidim (A) < \infty.$

\item Let $m$ be the minimum integer such that $A$ is a $m$-Gorenstein algebra, then: 

$$\fin (A) = \fidim (A) = \psidim (A)  = m.$$

\end{enumerate}

\end{teo}

The last result was obtained independently in \cite{kn:GS}.\\

An Artin  algebra $A$ is selfinjective if and only if $A$ is $0$-Gorenstein. Then as a consequence of the previous theorem we get the main result of \cite{kn:marchuard}. 

\begin{coro}\cite[Corollary 6]{kn:marchuard}
Let $A$ be an Artin  algebra. Then $A$ is selfinjective if and only if $\fidim(A) = 0$.
\end{coro}

In \cite{kn: LMM} has been proved that left and right $\phi$-dimensions for radical square zero algebras coincide but it is not the case for the $\psi$-dimension. As a consequence of Lemma \ref{Gorenstein reversible} and Corollary \ref{phi algGorenstein} we get the symmetry of $\phi$-dimension and $\psi$-dimension for Gorenstein algebras:

\begin{coro}
If $A$ is a  Gorenstein algebra, then: 
$$\fidim A = \fidim A^{op} = \psidim A = \psidim A^{op}$$
\end{coro}

We conjecture that the simetry of $\phi$-dimension holds for any Artin algebra.

\begin{ej}\label{EjGorenstein }

Let $A = \frac{\Bbbk Q}{I}$ be the following $\Bbbk$-algebra.

$$Q = \xymatrix{\dot 1 \ar@(l,u)^{\alpha} \ar[r]^{\gamma}&  2 \ar@(r,u)_{\beta}}$$\\

\noindent and $I = \langle \alpha^2, \beta^2, \gamma \alpha-\beta \gamma\rangle$.\\

Because $P(1) = I(2)$, $\id(P(2))= 1$ and $\pd (I(1)) = 1$, it follows that $A$ is a non selfinjective $1$-Gorenstein algebra. Therefore $\fidim A = \psidim A = 1$.

\end{ej}  

In the next subsections we compute the $\phi$-dimension and $\psi$-dimension for some widespread families of Gorenstein algebras.

\subsection{Gentle algebras}

A very important family of Gorenstein algebras are gentle algebras. 

\begin{defi}
Given a special biserial $\Bbbk$-algebra $A$, we call it \textbf{gentle algebra} if the following conditions are satisfied:

\begin{itemize}
  \item The ideal $I$ can be generated by paths of length $2$.
  \item For each arrow $\beta \in Q_1$ there exist at most one arrow $\gamma' \in Q_1$ such that $\beta \gamma' \in I$.
  \item For each arrow $\beta \in Q_1$ there exist at most one arrow $\alpha' \in Q_1$ such that $\alpha' \beta \in I$.
\end{itemize}
\end{defi}

In \cite{kn: GR} the authors proved that gentle algebras are $n$-Gorenstein, where $n$ is related to critical walks and gentle arrows. For definitions of critical walks and gentle arrows see \cite[Section 3.1]{kn: GR}.

\begin{teo}\cite[Theorem 3.4]{kn: GR}\label{Gentil->Gorenstein}
Let $A$ be a gentle algebra with $\eta(A)$ the maximum
length of critical walks starting with a gentle arrow, then

\begin{itemize}
  \item If $\eta(A) \geq 1$ then $\id_{A}A = \eta(A) = \pd_{A}D(A^{op})$, and $A$ is a $\eta(A)$-Gorenstein algebra.
  \item If $\eta(A) = 0$ then $\id_{A}A = \pd_{A}D(A^{op}) \leq 1$, and $A$ is selfinjective or a 1-Gorenstein algebra.
\end{itemize}

\end{teo}

\begin{coro}

 Let $A$ be a gentle algebra with $\eta(A)$ the maximum
length of critical walks starting with a gentle arrow, then

\begin{enumerate}
\item if $\eta(A) = 0$, $\fidim (A) = \psidim (A) \leq 1.$
\item if $\eta(A) \geq 1$, $\fidim (A) = \psidim (A) = \eta(A).$
\end{enumerate}

\end{coro}

\subsection{Cluster-tilted algebras}

An algebra $A$ is called cluster-tilted algebra if $A \cong \enn_{\mathcal{C}}(T)$, where $\mathcal{C}$ is a cluster category and $T$ is a tilting object in such category (see \cite{kn:5autores}).\\

Because each cluster-tilted algebra is $1$-Gorenstein (see \cite{kn:KR}), it follows the next corollary:

\begin{coro}

Let $A$ be a non selfinjective cluster-tilted algebra, then $\phi \dim (A) = \psi \dim (A) = 1$.

\end{coro}

The next example shows that there exist $1$-Gorenstein algebras which are not cluster-tilted.

\begin{ej}

Let $A$ be the path algebra $\frac{\Bbbk Q}{I}$, where $Q$ is the following quiver:

$$
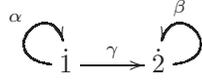
\xymatrix{\dot 1 \ar@(l,u)^{\alpha} \ar[r]^{\gamma}  & \dot 2 \ar@(r,u)}^{\beta}$$\\

\noindent and $I$ the ideal generated by $\{\alpha^{2}, \beta^{2}, \gamma \alpha- \beta \gamma\}$.\\

We get that $\fidim(A) = \psidim (A) = 1$ (see Example \ref{EjGorenstein }), but $A$ is not a cluster-tilted algebra because $\Omega^2\tau S_2 \ncong S_2$ (see \cite[Theorem 3.5]{kn:GS}). 

\end{ej}

\end{document}